\input amstex
\documentstyle{amsppt}
\magnification=\magstep1                        
\hsize6.5truein\vsize8.9truein                  
\NoRunningHeads
\loadeusm

\document

\topmatter

\title
Arestov's theorems on Bernstein's inequality 
\endtitle

\rightheadtext{Arestov's theorem on Bernstein's inequality in $L_p$ for all $p \geq 0}

\author Tam\'as Erd\'elyi
\endauthor

\address Department of Mathematics, Texas A\&M University,
College Station, Texas 77843, College Station, Texas 77843 \endaddress

\endthanks


\date April 26, 2019 
\enddate

\email terdelyi\@math.tamu.edu
\endemail

\abstract
We give a simple, elementary, and at least partially new proof of
Arestov's famous extension of Bernstein's inequality in $L_p$ to all $p \geq 0$.
Our crucial observation is that Boyd's approach to prove Mahler's inequality 
for algebraic polynomials $P_n \in {\Cal P}_n^c$ can be extended to all trigonometric 
polynomials $T_n \in {\Cal T}_n^c$.
\endabstract

\endtopmatter

\head 1. Introduction and Notation \endhead

Let ${\Cal T}_n^c$ be the collection of all trigonometric polynomials $T_n$ of the form
$$T_n(z) = \sum_{j=-n}^n{a_jz^j}\,, \qquad a_j \in {\Bbb C}\,, \quad z \in {\Bbb C} \setminus \{0\}\,.$$ 
Let ${\Cal P}_n^c$ be the collection of all algebraic polynomials $P_n$ of the form
$$P_n(z) = \sum_{j=0}^n{a_jz^j}\,, \qquad a_j \in {\Bbb C}\,, \quad z \in {\Bbb C}\,.$$
Let $D$ denote the open unit disk, and let $\partial D$ denote the unit circle 
of the complex plane. We define the Mahler measure (geometric mean of $Q$ on ${\partial D}$)
$$\|Q\|_0 = M_{0}(Q) := \exp\left(\frac{1}{2\pi} \int_{0}^{2\pi}{\log|Q(e^{it})|\,dt} \right)$$
for bounded measurable functions $Q$ on $\partial D$. It is well known, see [HL-52], for instance, that
$$\|Q\|_0 = M_0(Q) = \lim_{p \rightarrow 0+}{M_p(Q)}\,,$$
where
$$\|Q\|_p = M_p(Q) := \left( \frac{1}{2\pi} \int_{0}^{2\pi}{\left| Q(e^{it}) \right|^p\,dt} \right)^{1/p}\,, 
\qquad p > 0\,.$$
It is also well known that for a function $Q$ continuous on $\partial D$ we have
$$\|Q\|_\infty = M_{\infty}(Q) := \max_{t \in [0,2\pi]}{|Q(e^{it})|} = \lim_{p \rightarrow \infty}{M_p(Q)}\,.$$
It is a simple consequence of the Jensen formula that
$$M_0(Q) = |c| \prod_{k=1}^n{\max\{1,|z_k|\}}$$
for every polynomial of the form
$$Q(z) = c\prod_{k=1}^n{(z-z_k)}\,, \qquad c,z_k \in {\Bbb C}\,.$$

Bernstein's inequality  
$$\|T_n^\prime\|_\infty \leq n\|T_n\|_\infty\,, \qquad T_n \in {\Cal T}_n^c\,,$$
plays a crucial role in proving inverse theorems of approximation as well as many 
other results in approximation theory. See [BE-95], for instance. As far as the history of Bernstein's inequality 
is concerned we refer to Nevai's lovely papers [N-14] and [N-19] and the references in them. 
We do not repeat the full story here. In 1981 Arestov [A-81] proved that  
$$\|T_n\|_p \leq n\|T_n\|_p\,, \qquad T_n \in {\Cal T}_n\,, \tag 1.1$$ 
for all $p \geq 0$, extending the result known only for $p \geq 1$ for a long time. 
Simpler proof of Bernstein's inequality in $L_p$ for all $p \geq 0$ have been given by 
Golitschek and Lorentz in [GL-89] which is presented in the book [DL-93] by DeVore and Lorentz. 
A very elegant and even more simplified proof was published recently in [QZ-19] by 
Queff\'elec and Zarouf. A central part of their proof is to prove (1.1) for $p=0$ first.  
Mahler [M-61] showed (1.1) for $p=0$ but only for polynomials $P_n \in {\Cal P}_n^c$, and he gave 
a rather involved proof. Mahler's inequality was also posed as a problem by Vaaler in the Problems 
section of the American Mathematical Monthly and solved by Boyd [VB-91] using an elementary theorem 
of Bernstein. We note that Glazyrina [G-05] proved a sharp Markov-type inequality for algebraic 
polynomials in $L_0$ on finite subintervals of the real line.  

In this note we give a simple, elementary, and at least partially new proof of 
Arestov's famous extension of Bernstein's inequality in $L_p$ to all $p \geq 0$. 
Our crucial observation is that Boyd's approach to prove Mahler's inequality for 
algebraic polynomials $P_n \in {\Cal P}_n^c$ can be extended to all trigonometric polynomials 
$T_n \in {\Cal T}_n^c$.  

\proclaim{Theorem 1.1}
We have 
$$\|T_n^\prime\|_0 \leq n \|T_n\|_0\,, \qquad T_n \in {\Cal T}_n^c\,.$$  
Equivalently
$$\int_{\partial D}{\log|T_n^\prime(z)/n| \, |dz|} \leq \int_{\partial D}{\log|T_n(z)| \, |dz|}\,, 
\qquad T_n \in {\Cal T}_n^c\,.$$
\endproclaim

\proclaim{Theorem 1.2}
With the notation $\log^+|a| := \max\{\log|a|,0\}$ we have
$$\int_{\partial D}{\log^+|T_n^\prime(z)/n| \, |dz|} \leq \int_{\partial D}{\log^+|T_n(z)| \, |dz|}\,, 
\qquad T_n \in {\Cal T}_n^c\,.$$
\endproclaim

\proclaim{Theorem 1.3}
We have
$$\|T_n^\prime\|_p \leq n \|T_n\|_p\,, \qquad T_n \in {\Cal T}_n^c\,,$$
for every $p > 0$. 
\endproclaim

\head 2. Lemmas \endhead

To prove the Theorem 1.1 we need two lemmas.

\proclaim{Lemma 2.1} Associated with $S_n \in {\Cal T}_n^c$ let $P_{2n} \in {\Cal P}_{2n}^c$ be defined by 
$P_{2n}(z) = z^nS_n(z)$. If $P_{2n}$ has each of its $2n$ zeros in $D$,  then $S_n^{\prime}$ has each of its 
zeros in $D$ as well. The same is true if $D$ is replaced by the closed unit disk $\overline{D}$. 
\endproclaim

\demo{Proof}
We prove the lemma for $D$, the case of the closed unit disk $\overline{D}$ follows 
from this by a straightforward limiting argument.  
Suppose $a \notin D$, that is, $a \in {\Bbb C}$ and $|a| \geq 1$. 
Suppose also that $S_n \in {\Cal T}_n^c$ and $P_{2n}(z) := z^nS_n(z)$ has each of its $2n$ zeros in $D$, 
that is 
$$P_{2n}(z) = c \, \prod_{j=1}^{2n}{(z-z_j)}\,, \qquad 0 \neq c \in {\Bbb C}\,, \quad z_j \in D\,.$$
We have 
$$\frac{aS_n^{\prime}(a)}{S_n(a)} = \sum_{j=1}^{2n}{\frac{a}{a-z_j}} - \frac{an}{a} = 
\sum_{j=1}^{2n}{\frac{1}{1-z_j/a}} - n\,. \tag 2.1$$
Observe that $|z_j/a| < 1$ for each $j = 1,2,\ldots,2n$, and hence 
$$\text{\rm Re} \left( \frac{1}{1-z_j/a} \right) > \frac 12\,, \qquad j = 1,2,\ldots,2n\,. \tag 2.2$$
Combining (2.1) and (2.2) we obtain 
$$\text{\rm Re} \left( \frac{aS_n^\prime(a)}{S_n(a)} \right) > \frac{2n}{2} - n = 0\,.$$
We conclude that $S_n^\prime(a) \neq 0$. 
\qed \enddemo

\proclaim{Lemma 2.2} Associated with $V_n \in {\Cal T}_n^c$ let $R_{2n} \in {\Cal P}_{2n}^c$ be defined by
$R_{2n}(z) = z^nV_n(z)$ and suppose that $R_{2n}$ has each of its $2n$ zeros in the closed unit disk 
$\overline{D}$. If $T_n \in {\Cal T}_n$ and 
$$|T_n(z)| \leq |V_n(z)|\,, \qquad z \in \partial D\,, \tag 2.3$$
then
$$|T_n^\prime(z)| \leq |V_n^\prime(z)|, \qquad z \in \partial D\,.$$
\endproclaim

\demo{Proof} Without loss of generality we may assume that $R_{2n}$ has each of its $2n$ zeros in $D$, 
the case when $R_{2n}$ has each of its $2n$ zeros in the closed unit disk $\overline{D}$ follows 
from this by a straightforward limiting argument. Let $Q_{2n} \in {\Cal P}_{2n}^c$ be defined by 
$Q_{2n}(z) := z^nT_n(z)\,.$  Let $\alpha \in {\Bbb C}$, $|\alpha| < 1$, and  
$$S_n(z) := (V_n - \alpha T_n)(z) = z^{-n}(R_{2n} - \alpha Q_{2n})(z)\,. \tag 2.4$$ 
It follows from (2.3) that
$$|Q_{2n}(z)| \leq |R_{2n}(z)|\,, \qquad z \in \partial D\,,$$
and hence $|\alpha| < 1$ and the fact that $R_{2n} \in {\Cal P}_{2n}^c$ does not vanish on $\partial D$ 
imply that 
$$|\alpha Q_{2n}(z)| < |R_{2n}(z)|\,, \qquad z \in \partial D\,.$$
Therefore Rouche's Theorem implies that the polynomial $P_{2n} := R_{2n} - \alpha Q_{2n} \in {\Cal P}_{2n}^c$ 
and $R_{2n}$ has the same number of zeros in $D$, that is, $R_{2n} - \alpha Q_{2n}$ has each of its $2n$ zeros 
in $D$. By Lemma 2.1 and (2.4) we can deduce that $S_n^\prime = V_n^\prime - \alpha T_n^\prime$ 
has each of its zeros in $D$. In particular,
$$S_n^\prime(z) = V_n^\prime(z) - \alpha T_n^\prime(z) \neq 0, \qquad z \in \partial D\,,$$
for all $\alpha \in {\Bbb C}$, $|\alpha| < 1$. We conclude that 
$$|T_n^\prime(z)| \leq |V_n^\prime(z)|, \qquad z \in \partial D\,.$$ 
\qed \enddemo

\head 3. Proof of Theorems 1.1, 1.2, and 1.3 \endhead

\demo{Proof of Theorem 1.1}
Associated with $T_n \in {\Cal T}_n^c$ let $P_{2n} \in {\Cal P}_{2n}^c$ be defined by $P_{2n}(z) = z^nT_n(z)$. 
Without loss of generality we may assume that $P_{2n}$ has exactly $2n$ complex zeros, the case when 
$P_{2n}$ has less than $2n$ complex zeros follows from this by a straightforward limiting argument. 

\noindent Case 1. Suppose $P_{2n}$ has all its $2n$ zeros in the closed unit disk $\overline{D}$. 
It follows from Lemma 2.1 that $T_n^\prime$ has all its zeros in the closed unit disk $\overline D$. 
Let 
$$T_n(z) = \sum_{j=-n}^n{a_jz^j}\,, \qquad a_j \in {\Bbb C}\,, \quad z \in {\Bbb C} \setminus \{0\}\,.$$
Using Jensen's formula and the multiplicative property of the Mahler measure we can easily deduce that
$$\|T_n^\prime\|_0 = n \|T_n\|_0 = n|a_n|\,.$$

\noindent Case 2. Suppose that some of the zeros of $P_{2n}$ are outside the closed unit disk $\overline{D}$. 
Let $z_1,z_2,\ldots,z_m$ be the zeros of $T_n$ outside the closed unit disk $\overline{D}$ and 
let $z_{m+1},z_{m+2},\ldots,z_{2n}$ be the zeros of $T_n$ in the closed unit disk $\overline{D}$.    
We have
$$T_n(z) = a_nz^{-n}\prod_{j=1}^m{(z-z_j)} \prod_{j=m+1}^{2n}{(z-z_j)}\,.$$
We define 
$$V_n(z) := a_nz^{-n}\prod_{j=1}^m{(1-\overline{z_j}z)} \prod_{j=m+1}^{2n}{(z-z_j)}\,.$$ 
Observe that (2.3) holds, $V_n \in {\Cal T}_n^c$, and $R_{2n} \in {\Cal P}_{2n}^c$ defined by 
$R_{2n}(z) := z^nV_n(z)$ has each of its $2n$ zeros in the closed unit disk $\overline{D}$. 
Using Lemma 2.2, the (in)equality of the theorem in Case 1, and Jensen's formula, we obtain
$$\|T_n^\prime\|_0 \leq \|V_n^\prime\|_0 = n\|V_n\|_0 = n |a_n|\prod_{j=1}^m{|\overline{z_j}|} = 
n |a_n| \prod_{j=1}^m{|z_j|} = n\|T_n\|_0\,.$$    
\qed \enddemo

\demo{Proof of Theorem 1.2}
We follow the argument given in [QZ-19] to base our proof on Theorem 1.1. It is well-known, and by applying 
Jensen's formula it is ieasy to see, that 
$$\log^+|v| = \frac{1}{2\pi}\int_{\partial D}{\log|v+w| \, |dw|}\,, \qquad \quad v \in {\Bbb C}\,,$$
and hence
$$\log^+|v| = \frac{1}{2\pi} \int_{\partial D}{\log|v + wu| \, |dw|}\,, 
\qquad u \in {\partial D}\,, \quad v \in {\Bbb C}\,. \tag 3.1$$
Let $T_n \in {\Cal T}_n^c$, $w \in \partial D$, and $E_n(z) := z^n$.
Applying Theorem 1.1 with $T_n$ replaced by $T_n + wE_n$ we obtain
$$\int_{\partial D}{\log|T_n^\prime(z)/n + wE_{n-1}(z)| \, |dz|} \leq 
\int_{\partial D}{\log|T_n(z) + wE_n(z)| \, |dz|}\,.$$
Integrating both sides on $\partial D$ with respect to $|dw|$, then using Fubini's theorem and (3.1), 
we get the theorem
\qed \enddemo

\demo{Proof of Theorem 1.3}
We follow the argument given in [QZ-19] to base our proof on Theorem 1.2. Observe that
$$u^p = \int_0^\infty{\log^+(u/a)p^2a^{p-1} \, da}\, \qquad p > 0\,, \quad u \geq 0\,. \tag 3.2$$
Indeed, the integration by parts formula gives
$$\int_0^\infty{\log^+(u/a)p^2a^{p-1} \, da} = \int_0^u{\log^+(u/a)p^2a^{p-1} \, da}
= \int_0^u{pa^{p-1} \, da} = u^p\,.$$
For the sake of brevity we will use the notation $d\mu(a) = p^2a^{p-1} \, da$.
Using (3.2), Fubini's theorem, Theorem 1.2, and Fubini's theorem again, we obtain
$$\split \int_{\partial D}{\left| T_n^\prime(z)/n \right|^p \, |dz|} = & 
\int_{\partial D}{\left( \int_0^\infty{\log^+|T_n^\prime(z)/(na)| d\mu(a)} \right) \, |dz|} \cr
= & \int_0^\infty{\left( \int_{\partial D}{\log^+|T_n^\prime(z)/(na)| |dz|} \right) \, d\mu(a)} \cr
\leq & \int_0^\infty{\left( \int_{\partial D}{\log^+|T_n(z)/a| |dz|} \right) \, d\mu(a)} \cr
= & \int_{\partial D}{\left( \int_0^\infty{\log^+|T_n(z)/a| d\mu(a)}\right) \, |dz|} \cr
= & \int_{\partial D}{\left| T_n(z) \right|^p \, |dz|} \,. \cr \endsplit$$
\qed \enddemo

\head 4. Acknowledgement \endhead
The author thanks Herve Queff\'elec and Paul Nevai for checking the details 
of the proof in this paper and for their suggestions to make the paper more 
readable.

\medskip

\enddocument